\documentclass[12pt]{amsart}

\usepackage[centertags]{amsmath}
\usepackage{amsfonts}
\usepackage{amsthm}
\usepackage{newlfont}
\usepackage{amscd}
\usepackage{amsgen}
\usepackage{amssymb}
\usepackage{longtable}
\newlength{\defbaselineskip} 

\theoremstyle{plain}
\newtheorem{thm}{Theorem}[section]

\newtheorem{lem}{Lemma}[section]

\theoremstyle{remark} \newtheorem{rem}{Remark}[section]
\theoremstyle{definition} 
\theoremstyle{definition} \newtheorem{pro}{Problem}[section]
\theoremstyle{definition} \newtheorem{ex}{Example}[section] %

 \numberwithin{equation}{section}
\DeclareMathOperator{\Pic}{Pic}
 
 \DeclareMathOperator{\Def}{Def}
\DeclareMathOperator{\sing}{sing}

\begin{document}

\title{Primitive contractions of Calabi--Yau threefolds II}

\author{Grzegorz Kapustka}

\thanks{Research supported by the European Social Fund and the Polish government in
the frame of The Integrated Regional Operational Programme.} \maketitle
\begin{abstract} We construct several new examples of Calabi--Yau
threefolds with Picard group of rank $1$. Each of these examples
is obtained by smoothing the image of a primitive contraction with
exceptional divisor being a del Pezzo surface of degree $4$, $5$,
$6$, $7$ or $\mathbb{P}^1\times \mathbb{P}^1$.
\end{abstract}
\section{Introduction}
The lack of understanding of Calabi--Yau threefolds is the main
gap in the classification of $3$-dimensional algebraic varieties
(see \cite{Wilson1}, \cite{Wilson2}, \cite{Wilson}, \cite{namikawa
top}, \cite{Gross2}). For various reasons manifolds with Picard
group of rank 1 play a special role among all Calabi--Yau
threefolds (see \cite{Gross2}, \cite{EV}, \cite{Reid},
\cite{Lee1}, \cite{KN}, \cite{Lee2}, \cite{Bertin}). By the ``Reid
fantasy" conjecture they are expected to lead, via sequences of
conifold transitions, to all other Calabi--Yau threefolds.  The
already known examples of such threefolds are obtained as complete
intersections in homogeneous varieties,  by the method of Kawamata
and Namikawa (see \cite{KN}, \cite{Lee2}) based on smoothing
normal crossing varieties, or obtained by covering singular Fano
threefolds (see \cite{Lee1}, \cite{T}). The aim of this paper
 is to apply a method inspired by the ``Reid fantasy" to construct
 new families of Calabi--Yau threefolds with
 Picard group of rank $1$ that are not of any of the standard
 type. Four of them where predicted by van Enckevort and van
 Straten in \cite{EV}.

The method of construction can be described as follows. We first
find a singular Calabi--Yau threefold $X'$ that contains a del
Pezzo surface $D'$. Then consider the following schematic diagram:
\begin{center}
\setlength{\unitlength}{0.9mm}
\begin{picture}(60,30)(0,0)
\put(1,10){\vector(1,0){11}} \put(24,20){\vector(-1,-1){8}}
\put(28,20){\vector(1,-1){8}} \put(52,10){\vector(-1,0){12}}

\put(-2,10){\makebox(0,0){$\mathcal{M}$}}
\put(14,10){\makebox(0,0){X'}} \put(26,22){\makebox(0,0){X}}
\put(38,10){\makebox(0,0){Y}}
\put(54,10){\makebox(0,0){$\mathcal{Y}$}}
\end{picture}
\end{center}
where the morphism $X\rightarrow X'$ is an appropriate resolution
of singularities such that $X$ is a smooth Calabi--Yau threefold
and the strict transform $D$ of $D'$ is isomorphic to $D'$. The
morphism $X\rightarrow Y$ is a primitive contraction of type II
(see \cite{Wilson2}) having $D$ as exceptional divisor. The arrow
$\mathcal{M}\rightarrow X'$ (resp. $Y\leftarrow\mathcal{Y}$) means
that $X'$ can be smoothed inside the moduli space of smooth
Calabi--Yau threefolds $\mathcal{M}$ (resp.~$Y$ inside
$\mathcal{Y}$). The result of the construction is a generic
Calabi--Yau threefold from the family $\mathcal{Y}$.

Our constructions of Calabi--Yau threefolds $X'$ are analogous to
the ones in \cite[Sec.~5]{GM}. We first embed the del Pezzo
surface $D'$ by the anti-canonical embedding into a linear
subspace of a projective space (or a homogeneous space). Then,
choosing an appropriate hypersurface from the ideal of the del
Pezzo surface we construct a singular Calabi--Yau threefold $X'$
containing $D'$. Next, using the results of \cite{DH}, we prove
that the resulting Calabi--Yau varieties have only ODP
singularities. In fact Theorem \ref{nodal} shows in general that
an appropriate generic complete intersection has only ODP
singularities. The second problem is to prove that the strict
transform $D$ of the del Pezzo surface $D'\subset X'$ in the
resolution $X$ is the exceptional locus of a primitive contraction
of type II. This will be done if we prove that the rank of the
Picard group of $X$ is $2$. The last statement is proved in
Theorem \ref{rank Pic} using a generalization of the
Grothendieck--Lefschetz theorem given in \cite{RS}.

After smoothing the images of the constructed primitive
contractions, we obtain several examples of Calabi--Yau threefold.
From the results of \cite{GM} we compute the Hodge numbers of
$\mathcal{Y}_t$ (a generic element of $\mathcal{Y}$) and compare
them with the Hodge numbers of known families. We also find other
important invariants of the threefolds obtained: the degree of the
second Chern class, i.e.~$c_2\cdot H$, and the degree of the
generator of the Picard group, i.e.~$H^3$, where $H$ is the
generator of the Picard group of $\mathcal{Y}_t$. For each of the
first $13$ Calabi--Yau threefolds $\mathcal{Y}_t$ presented in
Table \ref{table1}, the linear system $|H|$ is very ample. In each
case we have $h^{1,1}(\mathcal{Y}_t)=1$.

\begin{center}
            \renewcommand*{\arraystretch}{1}
\begin{table}[h]\label{table1}
 ×

\caption{}
\begin{tabular}{|c|c|c|c|c|c|c|}

\hline

          No & $\deg D'$&$X'$&$\sing X'$&$\chi(\mathcal{Y}_t)$&$H^3$&  $h^0(H) $               \\ \hline
                                                                         \hline
       1& $4$&$X'_{2{,}2{,}2{,}2}\subset \mathbb{P}^7$                           &$12$ ODP&$-120$  &$20$&$9$                          \\ \hline
         2&   $5$&$X'_{2{,}2{,}2{,}2}\subset \mathbb{P}^7$                           &$18$ ODP&$-102$  &$21$&$9$                          \\ \hline
           3& $6$&$X'_{2{,}2{,}2{,}2}\subset \mathbb{P}^7$                           &$24$ ODP&$-84$   &$22$&$9$                              \\ \hline
4&$6$&$X'_{2{,}2{,}2{,}2}\subset \mathbb{P}^7$ &$24$ ODP&$-86$
&$22$&$9$                             \\ \hline
  5&$7$&$X'_{2{,}2{,}2{,}2}\subset \mathbb{P}^7$ &$30$
ODP&$-72$ &$23$&$9$                             \\ \hline

6&$4$&$X'_{2{,}2{,}3}\subset \mathbb{P}^6$&$16$ ODP&$-128$  &$2$
or $16$&  $8$ \\ \hline
 7&           $5$&$X'_{2{,}2{,}3}\subset \mathbb{P}^6$                               &$23$ ODP&$-108$  & $17$&$8$ \\ \hline
  8&          $6$&$X'_{2{,}2{,}3}\subset \mathbb{P}^6$                               &$30$ ODP&$-88$   & $18$&$8$                                  \\ \hline
 9&$6$&$X'_{2{,}2{,}3}\subset \mathbb{P}^6$                               &$30$ ODP&$-90$   & $18$&$8$                                 \\ \hline

   10&           $5$&$X'_{2{,}2{,}1}\subset G(2{,}5)$               &$15$ ODP&$-100$  &  $25$&$10$                                 \\ \hline
      11&      $5$&$X'_{1{,}1{,}3}\subset G(2{,}5)$               &$20$ ODP&$-120$  &  $20$&$9$                                 \\ \hline
 12&$6$&$X'_{1{,}1{,}1{,}1{,}2}\subset G(2{,}6)$        &$12$ ODP&$-96$  &   $34$&$12$                                \\ \hline
    13&     $6$&$X'_{1{,}1{,}1{,}1{,}2}\subset G(2{,}6)$        &$12$ ODP&$-98$  &     $34$&$12$                              \\ \hline

       14&      $8$&$X'_5\subset\mathbb{P}^4$                       &$24$ ODP&$-156$ &  $6$ or $48$&  sec. \ref{embed quadric}                           \\ \hline

          15&  $8$&$X'_{3{,}3}\subset\mathbb{P}^5$                 &$16$ ODP&$-116$
            &
           $10$ or $80$ &  sec. \ref{embed quadric}                               \\ \hline

            16&$8$&$X'_{2{,}2{,}3}\subset\mathbb{P}^6$             &$14$ ODP&$-120$ &    $13$ or $104$& sec. \ref{embed quadric}                             \\ \hline
            17&$8$&$X'_{2{,}2{,}2{,}2}\subset\mathbb{P}^7$         &$13$ ODP&$-106$ & $17$ or $136$&   sec. \ref{embed quadric}                    \\ \hline
           18& $1$&$X'_{3{,}6}\subset\mathbb{P}(1{,}1{,}1{,}2{,}3)$&$4$  ODP&$-256$ &\multicolumn{2}{|c|}{ $Y_{2,6}\subset \mathbb{P}(1,\dotsc,1,3)$ } \\ \hline
            19&$2$&$X'_{3{,}4}\subset\mathbb{P}(1{,}1{,}1,1,1{,}2)$&$8$  ODP&$-176$ &\multicolumn{2}{|c|}{ $Y_{2,4}\subset\mathbb{P}^5$}             \\ \hline
            20&$2$&$X'_6\subset\mathbb{P}(1{,}1,1,1{,}2)$          &$20$ ODP&$-200 $& \multicolumn{2}{|c|}{$Y_5\subset\mathbb{P}^4$}                 \\ \hline

    \end{tabular}
 \end{table}
            \end{center}
The Calabi--Yau threefolds number 2, 8, 10, 15 where predicted in
\cite{EV} and number 3, 4, 5, 9, 12, 13, 17 are new. The others
Calabi--Yau threefolds are discussed in Remarks \ref{2.5},
\ref{2.6}, \ref{2.7}, and \ref{2.8}. The values of $h^0(H)$ in
Table \ref{table1} where found with the help of the computer
algebra system Singular.

 The second aim of this paper is to determine which del Pezzo
surfaces can be the exceptional locus of a primitive contraction
of type II. From \cite{Gross1} we know that the exceptional locus
of such a contraction is either a normal Gorenstein rational del
Pezzo surface or a non normal Gorenstein surface (i.e.~such that
$\omega_E^{-1}$ is ample) of degree $7$, or a cone over an
elliptic curve of degree $\leq 3$. We show in Theorem \ref{result}
that each of the above surfaces except $\mathbb{P}^2$ blown up in
one point can be the exceptional locus of a primitive contraction
of type II.

\section*{Acknowledgements} This paper is part of my doctoral
thesis. I would like to express my gratitude to my advisor S.~Cynk
for his guidance. I thank M.~Kapustka for his enormous help. I
thank M.~Gross for suggesting me the method of embedding in
Section \ref{sec. 4 quadric}. I would like to thank
J.~Buczy\'{n}ski, A.~Langer, V.~Nikulin, S.~Pli\'{s}, P.~Pragacz,
M.~Reid, J.~Wi\'{s}niewski for answering questions, discussions,
and advice. I thank the referee for useful remarks.
 \section{Nodal threefolds}\label{section st6,7}
In this section we show a method that permits us to check if a
given complete intersection is nodal.

\begin{thm}\label{nodal} Let $D \subset \mathbb{P}^{s+2}$ be a smooth surface of codimension $s$ that is the scheme-theoretic
base locus of a linear system of hypersurfaces of degree $d$. Then
a generic complete intersection of $s-1$ hypersurfaces of degree
$d$ containing $D$ is a nodal threefold.
\end{thm}
\begin{proof}
First observe that a generic complete intersection $G$ of $s-2$
hypersurfaces $g_1,\dotsc,g_{s-2}$ of degree $d$ is smooth.
Indeed, it follows from the Bertini theorem that it is smooth away
from $D$. Singular points on $D$ appear when the rank of the
jacobian matrix of $g_1,\dotsc,g_{s-2}$ is smaller than $s-2$. As
the differentials $dg_1,\dotsc,dg_{s-2}$ correspond to global
sections of the globally generated vector bundle $\mathcal{I}_{S}
/\mathcal{I}^2_{S}(d)$, we conclude that the singular points
appear on $c_3(\mathcal{I}_{D} /\mathcal{I}^2_{D}(d))$ but this is
zero for dimensional reasons (see \cite[Example 14.4.3]{Fulton}).

Denote by $\pi\colon\tilde{G} \rightarrow G$ the blow up of $G$
along $D$, and by $E$ its exceptional divisor. Now $D\subset G$ is
a scheme-theoretic base locus of the linear system $\Lambda$ on
$G$ induced by the hypersurfaces of degree $d$ in
$\mathbb{P}^{s+2}$ containing $D$. It follows that the linear
system $\pi^*\Lambda-E$ is base-point-free. So the strict
transform of a general element of $\Lambda$ is nonsingular.

We claim that the singularities of a generic element $L$ of
$\Lambda$ are exactly the points above which the strict transform
of $\Lambda$ on $\tilde{G}$ contains a fiber of $\pi|_E$. Indeed,
it is enough to show the above fact in local analytic coordinates
around a point of $D$. To do this we follow \cite[Claim 2.2]{DH}.

 We claim that a generic element of $\pi^*\Lambda-E$ cuts $E$ along a
finite number of fibers of $\pi|_E$. Consider the morphism
$$\varphi\colon\ \tilde{G} \longrightarrow \mathbb{P}(H^0(\mathcal{O}_{\tilde{G}}(\pi^*\Lambda-E)))=:\mathbb{P}^N$$
given by $|\pi^*\Lambda-E|$. It follows from the assumptions that
the images of the fibers of $\pi|_E$ are lines in $\mathbb{P}^N$,
so the codimension of the set of divisors of the system $\Lambda$
that have a singularity at a given point of $D$ is $2$. We deduce
that a generic element of $\Lambda$ cannot have singularities
along a curve on $D$ (see \cite[p.~75]{DH}).

To finish the proof observe that the fibers of $\pi$ over $D$ are
isomorphic to $\mathbb{P}^1$. It follows that a generic element of
$\Lambda$ has a small resolution with $\mathbb{P}^1$ as
exceptional divisor.
 Now using the Bertini theorem on $E$ for the linear system
$(\pi^*\Lambda-E)|_E$, observe that a generic element of
$\pi^*\Lambda-E$ cuts $E$ along a nonsingular surface. It follows
that the normal bundle to an exceptional curve $C$ of the small
resolution has normal bundle with subbundle $\mathcal{O}_{C}(-1)$
(corresponding to the smooth strict transform of $D$).

From \cite[Thm.~4.4]{KOL} we know that a generic element $L$ has
only singularities of $cA$ type that are normal (hypersurface
singularities that are regular in codimension 1). In small
resolutions of such singularities there are only rational curves
with normal bundles $\mathcal{O}(-1)\oplus \mathcal{O}(-1)$ or
$\mathcal{O}\oplus \mathcal{O}(-2)$ (see \cite[Rem.~1.7]{Fr}).
Thus the normal bundle of an exceptional curve is
$\mathcal{O}(-1)\oplus \mathcal{O}(-1)$ in our cases, so the
singularities of a generic $\Lambda$ are ordinary double points
(see \cite{Reid 3}).

J.~Wi\'{s}niewski remarked that we can compute directly the normal
bundle of a contracted curve $C$. Indeed, the normal bundle
$N_{E|\tilde{G}}$ is $\mathcal{O}_E(-1)$ (since $\pi$ is the blow
up along $D$), and as observed above, $C$ is contained in a smooth
surface that is contained in $E$, which gives the other
summand.\end{proof}
\begin{rem}\label{chern nodes}
To compute the number of nodes on $L$ we proceed as follows. From
the Bertini theorem the double points of $L$ lie on $D$. More
precisely, since $L$ is a complete intersection of $s-1$
hypersurfaces, the corresponding sections of $\mathcal{I}_{D}
/\mathcal{I}^2_{D}(d)$ are linearly dependent precisely at the
singular points of $L$. This means that the number of nodes is
equal to $c_2(\mathcal{I}_{D} /\mathcal{I}^2_{D}(d))$.
\end{rem}
\begin{rem} It is interesting to compare the above theorem with
the Namikawa theorem: A Calabi--Yau threefold with terminal
singularities can be deformed to a Calabi--Yau threefold with
ordinary double points (see \cite{namikawa top}).
\end{rem}

\subsection{Embedding in a complete intersection of four
quadrics}\label{sec. 4 quadric} Let $D'_i \subset \mathbb{P}^i$
for $i=4,5,6,7$ be smooth del Pezzo surfaces of degree $i$
embedded by the anti-canonical system. We choose in $\mathbb{P}^7$
a linear subspace $L_i$ of dimension $i$ such that $D'_i \subset
L_i \subset \mathbb{P}^7$. Let $X'_i$ denote the intersection of
four generic quadrics from the ideal defining $D'_i$ in
$\mathbb{P}^7$. From Theorem \ref{nodal} the threefold $X'_i$ is a
nodal Calabi--Yau threefold.

Blowing up $D'_i$ we resolve the singularities of $X'_i$, and
flopping the exceptional divisors we obtain a smooth Calabi--Yau
threefold $X_i$ such that the strict transform $D_i$ of $D'_i$ is
isomorphic to $D_i'$. Denote by $H^*$ the pull-back to $X_i$ of
the hyperplane section $H$ of $\mathbb{P}^7$.
\begin{thm}\label{rank Pic} The rank of the Picard group of $X_i$ is $2$. Moreover, this group is generated by $D_i$ and $H^{*}$.
\end{thm}
\begin{proof} Observe first that $\Pic(X_i)$ is isomorphic in a natural way
to the Picard group of the blowing up $\tilde{X_i'}$ of $X_i'$
along $D_i'$.  Let $R_i$ be a generic complete intersection of
three quadrics containing $D'_i$. From the proof of Theorem
\ref{nodal} it follows that $R_i$ is smooth. Let
$\pi_i\colon\tilde{\mathbb{P}}^7\rightarrow \mathbb{P}^7$ be the
blow up along $D'_i$. The strict transforms of quadrics containing
$D'_i$ define a base-point-free linear system $\Lambda_i$ on the
strict transform $\tilde{R_i}$ of $R_i$. Observe that
$\tilde{X_i'}\in \Lambda_i$. We shall show that $\Lambda_i$ is
also big. Indeed, in the cases $i=4,5,6$ this is clear since
linear forms already separate points on $\mathbb{P}^7-L_i$.

\emph{Claim}. The system $|2H-D'_7| $ on $R_7-D'_7$ of quadrics
containing $D'_7$ is big and does not contract any divisor to a
curve.

Recall first that the ideal of $D_7\subset \mathbb{P}^7$ is
scheme-theoretically defined by the $2\times 2$ minors of a
$3\times 4$ matrix
$$\left( \begin{array}{ccc}
\begin{array}{ccc|c}
t&x&y&s\\
x&l&z&r\\
y&z&u&v\\
\end{array}

\end{array}\right)$$
obtained by deleting the last row from a symmetric matrix with
generic linear forms in $\mathbb{P}^7$. Since the first syzygies
of this set of quadrics are generated by linear forms we can use a
theorem of Room (see \cite[Prop.~3.1]{AR}). Since $(R_7 \cap
Sec(D'_7))\neq R_7\subset \mathbb{P}^7$ (where $Sec(.)$ denotes
the secant variety) we find that $|2H-D_7'|$ gives a birational
morphism. Moreover, the closure of an at least 2-dimensional fiber
of $|2H-D_7'|$ is a linear space $W$ that cuts $D'_7$ along a
curve of degree 2. It follows that the dimension of $W$ is 2 and
that $W$ cuts $D'_7$ along a curve from the system $|L-E_1|$
($D'_7$ is a blow up of $\mathbb{P}^2$ in two points, $L\subset
\mathbb{P}^2$ is a generic line, and $E_1$ and $E_2$ are
exceptional divisors). Suppose now that $|2H-D_7'|$ contracts a
divisor on $R_7$ to a curve. We obtain a contradiction since $R_7$
does not contain the three-dimensional scrolls determined by
$|L-E_1|$ and $|L-E_2|$.

We conclude that the linear system $\Lambda_i$ is big and defines
a morphism that contracts no divisor on $\tilde{R_i}$ to a point
or to a curve. From \cite[Thm.~6]{RS} one has an exact sequence
$$0\rightarrow K\rightarrow\Pic(\tilde{R_i})\rightarrow
\Pic(\tilde{X_i'})\rightarrow Q\rightarrow 0$$ where $K$ is the
subgroup generated by the divisors in $\tilde{R_i}$ which map to
points under the map given by $\Lambda$, and $Q$ is the group
generated by irreducible components of the traces on
$\tilde{X_i'}$ of those divisors on $\tilde{R_i}$ that map to a
curve under this map. We obtain $\Pic(\tilde{R_i})\simeq
\Pic(\tilde{X_i'})$. Now, arguing as before we conclude that
$\Pic(\tilde{R_i})$ is isomorphic to the Picard group of the
strict transform in $\tilde{\mathbb{P}^7}$ of an intersection of
two quadrics containing $D_i'$. So it is also isomorphic to
$\Pic(\tilde{\mathbb{P}}^7)$ which is generated by the exceptional
divisor and the pull-back of a hyperplane.
\end{proof}

\begin{rem}
We can also compute the number $h^{1,2}(X_i)$ in Singular
\cite{GPS} using the method described in \cite[Rem.~4.1]{GP}. From
the Euler characteristic of $X_i$, we can also deduce that
$h^{1,1}(X_i)=2$.
\end{rem}
Since $D_i|_{D_i}=K_{D_i}$ and $H^{*}|_{D_i}=K_{D_i}$ using
\cite[Lem.~2.5]{GM} we conclude that $D_i$ is the exceptional
divisor of a primitive contraction of type II.
\begin{thm}\label{thm Singular} The morphism $\varphi_{|H^{*}+D_i|}\colon X_i \longrightarrow
Y_i$ associated to the linear system $|H^{*}+D_i|$ is a primitive
contraction of type II with exceptional divisor $D_i$. Moreover,
$Y_i$ is smoothable and if $i=6$ there are two different
smoothings.
\end{thm}
\begin{proof}
We can show using \cite[Lem.~2.5]{GM} that some multiple
$n(H^{*}+D_i)$ gives a primitive contraction. To prove that it is
enough to take $n=1$, we shall show that $|H^{*}+D_i|$ is
base-point-free. Indeed, it is enough to find a divisor $G_i\in
|H^{*}+D_i|$ such that $G_i\cap D_i=\emptyset$. We first construct
a surface $G_i'\subset X_i'$. A generic quadric (or linear form)
that contains $D_i'\subset \mathbb{P}^7$ cuts $X_i'$ along a
divisor $D_i'+S_i$, where $S_i$ is an irreducible surface passing
through all the singularities of $X_i'$ (this follows from the
fact that $D_i'$ is smooth and there is no smooth Cartier divisor
passing through a singular point). With Singular \cite{GPS}, we
compute that there exists a cubic (resp.~a quadric) that contains
$S_i$ but does not contain $D_i'$. This cubic cuts $X_i'$ along
the divisor $S_i+G_i'$. Using Singular, we show that $G_i'\cap
D_i'$ is a finite set of points (these will be all the singular
points of $X_i')$. Moreover, we find that the ideal
$\mathcal{I}_{D'_i}+\mathcal{I}_{G'_i}$ is radical, so $G_i'$ and
$D_i'$ intersect transversally at these points. It follows that
their strict transforms $G_i$ and $D_i$ are disjoint.

We claim that $|H^{*}+D_i|$ gives a birational morphism that
contracts exactly $D_i$. Indeed, the system $|H^{*}|+D_i$ already
separates points from $X_i-D_i$ and $H^{*}+D_i|_{D_i}=0$. It
remains to prove the normality of the image $Y_i$.

Since $|H^*|$ separates points (also infinitely near) outside of
$D_i$, it is enough to prove that the image is normal at the point
$\varphi_{|G_i|}(D_i)$. We choose a generic hyperplane $K$ passing
through this point. It is enough to show that this hyperplane cuts
the image along a normal surface. Observe that the strict
transform of $K$ in $X_i$ is a smooth element  $H_1^* \in |H^*|$.
The divisor $ H_1^*$ is the strict transform of a hyperplane
section $H_1\subset X'_i$. Since $D'_i\cap H_1\subset H_1$ is a
linearly normal elliptic curve, the morphism
$$H^0(\mathcal{O}_{H_1\cap X'_1}(H)) \longrightarrow
H^0(\mathcal{O}_{D'_i\cap H_1}(H))$$ is a surjection. This implies
that
$$H^0(\mathcal{O}_{X_H}(G_i)\otimes \mathcal{O}_{X_H}(-D_i)) \longrightarrow
H^0(\mathcal{O}_Z(G_i)\otimes\mathcal{O}_Z
 (-D_i)),$$ where $Z=D_i\cap H_1^*$ and $X_H=X_i\cap H_1^*$, is a surjection.

To prove that $\varphi_{|Z|}(X_H)$ is normal, we can now argue as
in \cite[p.~416]{Har} using the fact that
$H^0(\mathcal{O}_Z(-D_i))=H^0(\mathcal{O}_Z(-Z))= m_{P'}/
m_{P'}^2$ (see \cite[Thm.~4.23]{Reid ch}) where $P'$ is the
elliptic singularity obtained after the contraction of $Z$.

To prove that $Y_i$ can be smoothed we use
\cite[Thm.~5.8]{Gross1}. To prove that in the case $i=6$ there are
two different smoothings we use \cite[Thm.~10]{namikawa1} or look
carefully at the proof of \cite[Thm.~4.3]{Gross1}. Indeed, the
versal Kuranishi space of a cone over a del Pezzo surface of
degree 6 has two components with different dimensions. It is
enough to observe that since $H^1_E (X, \Theta_{X})=0$ both these
components lead to smoothings of $Y_i$ (we use the local
cohomology sequence as in the proof of \cite[Thm.~3.3]{GM}).
\end{proof}

\begin{rem} The above proof gives an algorithm to compute
the number of nodes on $X_i'$; this number is equal to the degree
of $ (G'_i \cap D'_i)$ providing the latter is finite.
\end{rem}

 Using the results of \cite{GM} we can compute the Hodge numbers of
 smooth Calabi--Yau threefolds $\mathcal{Y}^i_t$ from smoothing families of $Y_i$. We obtain
 $h^{1,1}(\mathcal{Y}^i_t)=1$ and $h^{1,2}(\mathcal{Y}^i_t)=h^{1,1}(\mathcal{Y}^i_t)-
 \frac{1}{2}\chi(\mathcal{Y}^i_t)$. Denote by $T$ the generator of the
 Picard group of $\mathcal{Y}^i_t$. We compute the degree $T^3$ using the fact that $(H^{*}+D_i)^3$
 is not divisible by a cube and that the image of $H^{*}+D_i$ in $Y^i$ is a Cartier
 divisor (see the proof of \cite[Prop.~3.1]{GM}).

\begin{center}
            \renewcommand*{\arraystretch}{1}
\begin{table}[htp]
 ×

\caption{}
\begin{tabular}{|c|c|c|c|c|}

\hline
            $i$&The number of nodes on $X_i'$&$\chi(X_i)$&$\chi(\mathcal{Y}^i_t)$&$T^3=\deg(\mathcal{Y}^i_t)$                       \\ \hline
                                                                   \hline

            $4$&12&$-104 $&$-120$& $20$           \\ \hline
            $5$&18 &$-92$&$-102$& $21$ \\ \hline
            $6$& $24$&$-80$&$-84$ and $-86$&$22$                    \\ \hline
            $7$&$30$&$-64$&$-68$&$23$
            \\ \hline
            \end{tabular}
 \end{table}
            \end{center}
\begin{rem}\label{2.5} In the first case the obtained Calabi--Yau is the
intersection of the Grassmaniann $G(2,5)\subset \mathbb{P}^9$ with
a linear space and two quadrics. In the second case we obtain the
Calabi-Yau threefold of degree $21$ predicted in \cite{EV}.
\end{rem}

 It is also natural to find $\dim(|T|)$ or equivalently
$ \frac{T\cdot c_2}{12} =\chi(\mathcal{O}(T))-\frac{T^3}{6}$ in
the examples above. From the following two lemmas we deduce that
$\dim(|T|)=9$ in each case.
\begin{lem}\label{lem2.1} The dimension of the linear system $|G_i|$ is $9$.\end{lem}

\begin{proof} We have $G_i|_{G_i}=K_{G_i}$ and the image of $G_i$ in $Y_i$
is a hyperplane section of $Y_i$ (since
$h^1(\mathcal{O}_X(G))=0$). It follows that it is enough to prove
that $G_i'\subset \mathbb{P}^7$ (with the notation from the proof
above) is an embedding given by the canonical divisor. Since by
the previous computer calculations, $G_i\in |H^{*}+D_i|$, we have
that $K_{G_i}=(H^*+D_i)|_{G_i}=H^*|_{G_i}$ so
$H|_{G_i'}=K_{G_i'}$. Observe now that the system of cubics that
contains $S_i'$ and does not contain $D_i'$ is base-point-free on
$X_i'$ outside the singularities of $X_i'$. Indeed, we can add to
such a cubic a cubic obtained by multiplying a linear form with
the generator of the ideal of $D_i'\subset \mathbb{P}^7$ that
defines $S_i'$. It follows that a generic surface
$G_i'\subset\mathbb{P}^7$ constructed above is smooth outside the
singularities of $X_i'$. Since the ideal
$\mathcal{I}_{D'_i}+\mathcal{I}_{G'_i}$ is radical, $G_i'$ is
smooth everywhere. Finally, since the divisor $H^*-G$ is not
effective, the surface $G_i'\subset\mathbb{P}^7$ is not contained
in a hyperplane. It remains to show that $G_i'\subset\mathbb{P}^7$
is is embedded by a complete linear system. First, the surface
$D_i'\subset\mathbb{P}^7$ is arithmetically Cohen-Macaulay. Since
the last property is preserved via Gorenstein linkage (the general
properties of linkages are described in \cite{KM}), we obtain
$h^1(\mathcal{I}_{G_i'}(1))=0$. So $G_i'\subset\mathbb{P}^7$ is
linearly normal.
\end{proof}
Thus $Y_i\subset \mathbb{P}^8$ is linearly normal, so
$h^1(\mathcal{I}_{Y_i}(1))=0$. It remains to use the
semi-continuity theorem \cite[III 12.8]{Har} and the following
lemma.
\begin{lem} Let $i\colon Y\hookrightarrow \mathbb{P}^{n}$ be a
smoothable (normal Gorenstein) Calabi--Yau threefold. Then $Y$ can
be smoothed inside $\mathbb{P}^{n}$.
\end{lem}
\begin{proof}
 Let $\Def(Y|\mathbb{P}^{n})$ be the functor from the category of
germs of complex spaces to the category of sets, such that the
image of a germ $(S,x)$ is the set of isomorphism classes of
deformations
\[ \begin{CD}
             {Y}   @>i>>    {\mathbb{P}^{i+1}} @>>> {x}\\
             @VVV         @VVV     @VVV\\
             {\mathcal{Y}}   @>>> {\mathbb{P}^{i+1}\times
             S}@>>> {S}
             \end{CD}\]
of the embedding $Y\hookrightarrow \mathbb{P}^{n}$ (here
$\mathcal{Y}\rightarrow S$ is flat). We need to show that for each
germ $S$ the natural map
$$\Def(Y|\mathbb{P}^{n})(S)\rightarrow \Def(Y)(S)$$ is surjective.
The latter follows from Theorem 1.7 (ii) in \cite{Weh} if we show
that $H^1(\Theta_Y|_{\mathbb{P}^{n}})=0$.

 Indeed, since $\dim(\sing(Y))=0$  and $Y$ is normal Gorenstein it follows from
\cite[Prop.~1.1]{AJ} that $H^1(\mathcal{O}_Y(1))=0$. The claim
follows from the long exact cohomology sequence associated to
\begin{equation*}0\rightarrow \mathcal{O}_Y \rightarrow (\mathcal{O}_Y
(1))^{n+1}\rightarrow\Theta_{\mathbb{P}^{n}}|_Y\rightarrow 0.
\qedhere
\end{equation*}
\end{proof}
\begin{ex} Let us show how to compute the number of nodes in the
case $i=7$. From Remark \ref{chern nodes} we need to find
$c_2(\mathcal{I}_{D_7} /\mathcal{I}^2_{D_7}(2))$. We use the
exact sequence
 $$0\rightarrow \mathcal{T}_{D_7}  \rightarrow \mathcal{T}_{\mathbb{P}^7}|_{D_7}
  \rightarrow \mathcal{N}_{\mathbb{P}^7|D_7}
  \rightarrow 0 .$$
   Observe that $c_t(D_7)=1-ht+\chi_{\mathrm{top}}(D_7)Pt^2=1-ht+\frac{5}{7}h^2t^2$ and $c_t(\mathcal{T}_{\mathbb{P}^7}|_{D_7})=(1+ht)^8$ where $h$ is a
  hyperplane section in $\mathbb{P}^7$ and $P\in D_7$ a point.
  Consequently, $c_t(\mathcal{T}_{\mathbb{P}^7}|_{D_7})c_t^{-1}(D_7)=c_t(
  \mathcal{N}_{\mathbb{P}^7|D_7})$, and it remains to apply \cite[Rem.~3.2.3]{Fulton}.

\end{ex}

\begin{ex}

In the same way we can embed del Pezzo surfaces $D'$ of degree
$4$, $5$, and $6$ into nodal complete intersections of a cubic and
two quadrics (we first embed into a smooth complete intersection
of two quadrics then use Theorem \ref{nodal}). The Calabi--Yau
threefolds $\mathcal{Y}_t$ from the resulting smoothing families
have $h^{1,1}(\mathcal{Y}_t)=1$.
\begin{center}
            \renewcommand*{\arraystretch}{1}
\begin{table}[htp]
 ×

\caption{}
            \renewcommand*{\arraystretch}{1}
            \begin{tabular}{|c|c|c|c|c|}\hline
            $\deg(D')$&The number of nodes on $X'$&$\chi(X)$&$\chi(\mathcal{Y}_t)$&$\deg(\mathcal{Y}_t)$                       \\ \hline

            $4$&$16$ &$-112$&$-128$& $16$ or $2$  \\ \hline

             $5$&$23$ &$-98$&$-108$& $17$ \\ \hline

             $6$&$30$ &$-84$&$-88$ and $-90$& $18$ \\ \hline
            \end{tabular}
        \end{table}    
            \end{center}

 \begin{rem}\label{2.6} In the case when $\deg(D')=4$, the invariants are similar to the invariants of a complete
intersection of four quadrics, when $\deg(D')=5$ to one
constructed in \cite{Lee1}. The existence of a Calabi--Yau
threefold of degree $18$ with Euler characteristic $-88$ where
predicted in \cite{EV}.\end{rem}
\end{ex}
\subsection{Embedding in complete intersections in Grassmannians} We
can use a similar construction to find primitive embeddings of del
Pezzo surfaces into Calabi--Yau threefolds that are complete
intersections in Grassmannians.

Let $D_1'\subset  \mathbb{P}^9$ be a del Pezzo surface of degree
$5$ anti-canonically embedded in a codimension $4$ linear section.
It is known that $D_1'$ can be seen as a section of the
Grassmannian $G(2,5)$ embedded via the Pl\"{u}cker embedding.

Let $X_1'$ be the intersection of $G(2,5)$ with a hyperplane $L$
and two quadrics containing $D_1'$.
\begin{lem}
The threefold $X_1'$ is a nodal Calabi--Yau threefold that can be
resolved by blowing up $D_1'$. Flopping the exceptional curves of
this blow up, we obtain a small resolution $X_1\rightarrow X_1'$.
The rank of the Picard group of $X_1$ is $2$.
\end{lem}
\begin{proof} To prove that $X_1'\subset G(2,5) \cap L$ is nodal,
we follow the proof of Theorem \ref{nodal}. Since the system on $
G(2,5) \cap L$ of quadrics that contain $D_1'$ separates points on
$(G(2,5) \cap L)-D_1'$, we can argue as in the proof of Theorem
\ref{rank Pic} to obtain the second part. In particular we
construct the surface $G_1'$.
\end{proof}
 We compute that
 $X'$ has $15$ ordinary double points (the radical of the ideal $\mathcal{I}_{D_1'}+
 \mathcal{I}_{G_1'}$ is $0$-dimensional and has degree $15$).
The strict transform of $D_1'$ in $X$ is an exceptional locus of a
primitive contraction of type II.

 In the same way we can embed a del Pezzo surface $D_2'\subset \mathbb{P}^9$ of degree
 $5$ into a nodal Calabi--Yau threefold $X_2'$ with $20$ nodes that is a complete
 intersection of $G(2,5)$ with two hyperplanes and a cubic.
Let $X_2\rightarrow X_2'$ be the appropriate small resolution and
$D_2$ the strict transform of $D_2'$.
 We obtain as before a primitive contraction given by $H^*+D$, with $D_2$ as exceptional
 locus.

The del Pezzo surface $D_3'\subset \mathbb{P}^6$ of degree 6 can
be seen as a special linear section of the Grassmannian
$G(2,6)\subset \mathbb{P}^{14}$. More precisely, it can be
described by the $4\times 4$ Pfaffians of a $6\times 6$
extra-symmetric matrix (i.e.~skew-symmetric and symmetric with
respect to the other diagonal). We embed $D'_3$ into a nodal
Calabi--Yau threefold $X_3'\subset \mathbb{P}^{10}$ that is a
complete intersection of $G(2,6)$ with a quadric and four linear
forms. We compute that $X_3'$ has $12$ nodes.

 Denote by $Y_i$ the images of the above primitive contractions and
by
 $\mathcal{Y}^i_t$ their smoothing families. In each case, we have
 $h^{1,1}(\mathcal{Y}^i_t)=1$.
\begin{center}

            \renewcommand*{\arraystretch}{1}
\begin{table}[htp]
 ×

\caption{}
            \renewcommand*{\arraystretch}{1}
            \begin{tabular}{|c|c|c|c|c|}\hline
            $i$&The number of nodes on $X_i'$&$\chi(X_i')$&$\chi(\mathcal{Y}^i_t)$&$deg(\mathcal{Y}^i_t)$                       \\ \hline

           $1$&$15$&$-105$ &-100& $25$           \\ \hline

            $2$&$20$ &$-130$&$-120$& $20$ \\ \hline
            $3$&$12$&$-104$&$-96$ and $-98$& $34$ \\ \hline

            \end{tabular}
\end{table}            
            \end{center}

\begin{rem} \label{2.7}These computations suggest that the smooth elements
$\mathcal{Y}_t$ from the smoothing family of $Y_2$ are generic
complete intersections of $G(2,5)$ with a hyperplane and two
quadrics. However, calculations in Singular show that the surface
$G'$ (see the proof of Theorem \ref{rank Pic}), which is
isomorphic to a hyperplane section of $Y_2$, is not generated by
quadrics. The elements from the smoothing family of $Y_1$ have the
same invariants as a Calabi--Yau threefold constructed in
\cite{Lee1} (see also \cite{EV}).
\end{rem}

\subsection{Embeddings of quadrics }\label{embed quadric}
 We shall embed as before a quadric $D'\subset \mathbb{P}^3$
(this embedding is given by half of the canonical divisor on $D'$)
into a nodal quintic in $\mathbb{P}^4$, a nodal complete
intersection of two cubics in $ \mathbb{P}^5$, a complete
intersection of a cubic and two quadrics, and a complete
intersection of four quadrics.

Let $ X'\subset \mathbb{P}^5$ be a generic complete intersection
of two cubics containing $D'\subset L\subset \mathbb{P}^5$ where
$L$ is a codimension $2$ linear section. We check that the
assumptions of Theorem \ref{nodal} hold, so $X'$ is a nodal
Calabi--Yau threefold. The blow up of $D'\subset X'$ is a small
resolution. Denote by $X$ the flop of the exceptional curves of
this resolution and by $D\subset X$ the strict transform of $D'$.
We prove as in Theorem \ref{rank Pic} that the rank of the Picard
group of $X$ is $2$. It follows that $D$ is the exceptional locus
of a primitive contraction of type II.

We can perform the same construction for $X'\subset\mathbb{P}^N$
being a nodal quintic (resp. a complete intersection of four
quadrics and a complete intersection of a cubic and two quadrics).
 If $Y$ is the image of the constructed
contraction and $\mathcal{Y}_t$ the smoothing family of $Y$ we
obtain as before $h^{1,1}(\mathcal{Y}_t)=1$. In these cases
however the primitive contractions are given by the linear systems
$|2H^{*}+D|=|G|$ (because $2H^*|_D=-K_D$). To compute the
dimension $h^0(\mathcal{O}_X(G))$ we observe that $G|_G=K_G\simeq
2H|_{G'}$. Arguing as in the proof of Lemma \ref{lem2.1}, we
obtain $H^1(\mathcal{I}_{G'}(2))=0$. Now, the elements of degree
$2$ in the ideal $\mathcal{I}_{G'}$ are exactly the quadrics
containing $X'$ (use Singular). Since $h^0(\mathcal{I}_{G'}(2))$
is the dimension of the kernel of the map
$H^0(\mathcal{O}_{\mathbb{P^N}}(2H))\rightarrow
H^0(\mathcal{O}_{G'}(2H))$, we are ready to compute
$h^0(\mathcal{O}_X(G))$ in Table $5$.

\begin{rem}\label{2.8}
The Calabi--Yau threefold $\mathcal{Y}_t$ obtained in the first
row of Table $5$ has similar invariants as the complete
intersection $Y_{3,4}\subset \mathbb{P}(1,1,1,1,1,2)$. In the
second row we probably obtain the Calabi-Yau threefold of degree
$10$ predicted in \cite{EV} and in the third row the
 Calabi--Yau threefold defined by the $4\times 4$ Pfaffians
of a $5\times 5$ skew symmetric matrix that has one row and one
column with
 quadric entries.
\end{rem}
            \renewcommand*{\arraystretch}{1.3}
\begin{table}[htp]\label{table5}
 ×

\caption{}
            \renewcommand*{\arraystretch}{1}
            \begin{tabular}{|c|c|c|c|c|c|}\hline
            $\deg(X')$&ODP on $X'$&$\chi(X')$&$\chi(\mathcal{Y}_t)$&$(2H^*+D)^3$ &$h^0(2H^*+D)$                      \\ \hline

           $ 5$&$24$&$-200$ &$-156$& $48$ & $16$          \\ \hline

            $9$&$16$ &$-144$&$-116$&  $80$&$22$ \\ \hline
             $12$&  $14 $ &$-144$&$-120$& $104$&$27$                         \\ \hline
            $16$&$13$&$-128$&$-106$& $136$&$33$                              \\ \hline

            \end{tabular}
\end{table}

\subsection{Embeddings in weighted projective spaces}\label{weighted}
Let us show one more application of Theorem \ref{nodal} and
\ref{rank Pic}.
  Let $D'\subset \mathbb{P}(1,1,1,1,2,3)$ (with
coordinates $x,y,z,t,u,v$) be a del Pezzo surface of degree $1$
defined by $x=y=s=0$, where $s$ is a generic sextic in
$\mathbb{P}(1,1,1,1,2,3)$. Let $q=xf+yg$ where $f, g$ are generic
quadrics. Denote by $X'=X'_{3,6}$ the variety defined in
$\mathbb{P}(1,1,1,1,2,3)$ by the equations $s=q=0$.
\begin{lem} The threefold $X'$ is a nodal Calabi--Yau threefold with four nodes.
\end{lem}
\begin{proof} Consider the sixfold ramified covering
$$c\colon\mathbb{P}^5 \ni(x,y,z,t,u_1,v_1) \rightarrow (x,y,z,t,u_1^2,v^3_1)\in \mathbb{P}(1,1,1,1,2,3).$$
 Since a generic sextic omits the singularities of $\mathbb{P}(1,1,1,1,2,3)$, it follows from the Bertini theorem that
 the singularities of $X'$ can occur only when $x=y=0$. To prove that these are ordinary
double points it is enough to show that $c^{-1}(X')\subset
\mathbb{P}^5$ has only nodes. The threefold $c^{-1}(X')$ is
defined by the sextic $s_1=s(x,y,z,t,u_1^2,v_1^3)$ and a cubic of
the form
$$xq_1(x,y,z,t)+xu_1^2+yq_2(x,y,z,t)+yu_1^2.$$
Since the system of such cubics defines the surface $c^{-1}(D')$
scheme-theoretically on the fourfold $s_1=0$, we can conclude as
in the proof of Theorem \ref{nodal}.

To compute the number of nodes we find the degree of the Jacobian
ideal of $X'$ (using Singular). We deduce that $c^{-1}(X')$ has 24
nodes at general points of $D'$.
\end{proof}
Let $X$ be as before the small resolution of $X'$ such that the
strict transform $D$ of $D'$ is isomorphic to $D'$. Since the
system of cubics of the form $$xq_1(x,y,z,t)+xu+yq_2(x,y,z,t)+yu$$
induces a system that separates points on $W- D'$ we follow the
proof of Theorem \ref{rank Pic} to obtain $\rho(X)=2$. Hence we
can find as before a primitive contraction with exceptional locus
$D$.

We find analogous constructions for a del Pezzo surface of degree
$2$ embedded in a natural way into the nodal Calabi--Yau
threefolds $X'_{3,4}\subset \mathbb{P}(1,1,1,1,1,2)$ and
$X'_6\subset\mathbb{P}(1,1,1,1,2)$. Assuming $\rho(X')=2$, we
deduce, with the same notation as in the previous section, that
$h^{1,1}(\mathcal{Y}_t)=1$.
\begin{center}
            \renewcommand*{\arraystretch}{1.3}
\begin{table}[htp]
 ×
\caption{}
            \renewcommand*{\arraystretch}{1}
            \begin{tabular}{|c|c|c|c|c|}\hline
            $X'$&The number of nodes on $X'$&$\chi(X')$&$\chi(\mathcal{Y}_t)$&$\deg(\mathcal{Y}_t)$                       \\ \hline

            $X'_{3,6}$&$4$&$-200$&$-256$&$4$  \\ \hline

              $X'_{3,4}$&$8$&$-148$&$-176$&$8$ or $1$  \\ \hline

           $X'_{6}$&$20$&$-184$&$-200$&$5$  \\ \hline
            \end{tabular}
      \end{table}      
            \end{center}
 Denote by $G'$ the smooth surface in $X'$ defined
by $f=g=0$ and by $G$ its strict transform in $X$. The above
calculations suggest the following theorem.
\begin{thm}\label{weigth}
The linear system $|6G|$ gives a primitive contraction with $D$ as
exceptional locus into a Calabi--Yau threefold that is isomorphic
to a complete intersection $X_{2,6}\subset
\mathbb{P}(1,1,1,1,1,3)$.
\end{thm}
\begin{proof} Observe that $G'$ and $D'$ intersect transversally
at 24 general points defined by
$g(x,y,z,t,u_1^2)=f(x,y,z,t,u_1^2)=s_1=0$.

Moreover, $G'\in |H'+D'|$ (where $H'$ is an element from the
system $\mathcal{O}_{\mathbb{P}(1,1,1,1,2,3)}(1)$). It follows
that $G\cap D=\emptyset$. Since the system $|6H'|$ is very ample,
we can argue as in \cite[Prop.~4.3]{GM} to show that $|6G|$ gives
a birational morphism $\varphi_{|6G|}$. Observe that $G'$ is a
surface of general type naturally isomorphic to $X_{2,6}\subset
\mathbb{P}(1,1,1,1,3)$. It follows that $6G|_G$ gives a morphism
into a projectively normal surface. We deduce as in the proof of
\cite[Prop.~4.3]{GM} that the image of $\varphi_{|6G|}$ is
projectively normal with one singular point. Now, as in the proof
of \cite[Thm.~4.5]{GM}, we can compute that the Hilbert series of
$\bigoplus_{n=0}^{\infty} H^0(\mathcal{O}_X(nG))$ is equal to
$$P(t)=\frac{(1-t^2)(1-t^6)}{(1-t)^5(1-t^3)}.$$ We conclude that
the image of $\varphi_{|6G|}$ is isomorphic to a normal complete
intersection $X_{2,6}\subset \mathbb{P}(1,1,1,1,1,3)$ with an
isolated rational Gorenstein singularity.
\end{proof}
\begin{rem} An analogous theorem holds for del Pezzo surfaces od degree $2$ naturally embedded into $X'_{3,4}\subset
\mathbb{P}(1,1,1,1,1,2)$ or $X'_6\subset\mathbb{P}(1,1,1,1,2)$.
The image of the primitive contraction is then isomorphic to
$X'_{2,4}\subset \mathbb{P}^5$ or to $X'_5\subset \mathbb{P}^4$
respectively.
\end{rem}
\begin{rem}
Note that the we cannot give an analogous construction for del
Pezzo $D'\subset X'_{4,6} \subset \mathbb{P}(1,1,1,2,2,3)$ of
degree 1. In fact, $X'_{4,6}$ will have a curve of singularities.
\end{rem}
\subsection{Primitive contractions of singular del Pezzo
surfaces}

Let $D'\subset \mathbb{P}^3$ be a cubic surface with du Val
singularities (rational Gorenstein del Pezzo surface of degree
$3$).

Let $X'$ be the quintic defined in $\mathbb{P}^4$ by the equation
$r=cq+lp$ where $c$ is a cubic, $l$ a linear form such that
$c=l=0$ define $D'\subset \mathbb{P}^3$, $p$ is a generic quartic,
and $q$ a generic quadric.

\begin{lem}\label{lemma}
The threefold $X'$ is a nodal Calabi--Yau threefold with $24$
nodes at points lying outside the singularities of $D'$.
\end{lem}
\begin{proof} The partial derivatives of $r$ give a system of
equations
$$q\frac{dc}{dx_1}+c\frac{dq}{dx_1}+l\frac{dp}{dx_1}+p\frac{dl}{dx_1}=0$$
$$\dotsc $$
$$q\frac{dc}{dx_5}+c\frac{dq}{dx_5}+l\frac{dp}{dx_5}+p\frac{dl}{dx_5}=0.$$
From the Bertini theorem the singularities of $X'$ lie on $D'$.
For a generic choice of $p$ and $q$ the solution of this system is
the set of $24$ points where $p=q=l=c=0$. Indeed, we can assume
that $l=x_1$; then if we put $c=l=0$ the first equation takes the
form $p+q\frac{dc}{dx_1}=0$. We infer that the singularities of
$X'$ lie outside those of $D'$. Furthermore at points $P\in D'$
where $D'$ is nonsingular, one of the partial derivatives of the
equation defining $D'$, say $\frac{dc}{dx_2}$, is nonzero. Taking
the second equation we obtain $q(P)=0$, thus $p(P)=0$.
\end{proof}

Let $X$ be the Calabi--Yau threefold obtained by blowing up $X'$
along the smooth surface defined by $l=q=0$. It follows from
Theorem \ref{rank Pic} that the rank of the Picard group of $X$ is
2. The strict transform $D\simeq D'$ of $D'$ in $X$ is the
exceptional locus of a primitive contraction of type II.

Using the constructions from Sections \ref{embed quadric},
\ref{sec. 4 quadric}, and \ref{weighted}, we can find, arguing as
in \ref{lemma}, primitive contractions with exceptional loci being
rational Gorenstein del Pezzo surfaces of degree $1$ to $7$, and a
(possibly singular) quadric. The exceptional divisor of a
primitive contraction of type II is either a normal Gorenstein del
Pezzo surface with rational double points, or a nonnormal del
Pezzo surface of degree $7$ whose normalization is a nonsingular
rational ruled surface, or a Gorenstein del Pezzo surface of
degree $\leq 3$ that is a cone over an elliptic curve (see
\cite[Lem.~3.2]{Wilson}). The converse also holds.
\begin{thm}\label{result} If $E$ is a del Pezzo surface as above that is not isomorphic to the blow up of $\mathbb{P}^2$ in one point,
then we can find a primitive contraction of type II with
exceptional locus isomorphic to $E$.
\end{thm}
\begin{proof} It suffices to consider the case when $E$ is non-normal. Then from \cite{R1}
and \cite{Wilson}, $E$ is isomorphic to one of the following:

\begin{enumerate}
\item{} $\tilde{F}_{5;1}$ the projection from a point that lie in
the plane spanned by the line $l$ and a fiber of the rational
normal scroll that is a join in $\mathbb{P}^8$ of $l$ and a
rational normal curve of degree $6$, or \item{} $\tilde{F}_{3;2}$
the projection from a point on the plane spanned by the conic of
the scroll that is the join in $\mathbb{P}^8$ of a conic and a
rational normal curve of degree $5$.
\end{enumerate}
 We find that $\tilde{F}_{3;2} \subset \mathbb{P}^7$ is defined by
 $14$ quadrics. After a suitable change of coordinates $x,z,t,u,v,w,p,q$ they can be
 written as $p^2-wq
,wp-vq ,vp-uq ,up-tq ,w^2-uq ,vw-tq ,uw-tp ,zw-xq ,v^2-tp ,uv-tw
,zv-xp ,u^2-tv ,zu-xw ,zt-xv$. These quadrics define a big system
that does not contract any divisor to a curve on $R_i
-\tilde{F}_{3;2}$ ($R_i$ is a smooth intersection of three such
quadrics) so we can argue as before. We first choose four generic
quadrics from the ideal of $\tilde{F}_{3;2} \subset \mathbb{P}^7$,
and we denote by $X'$ their intersection. We prove as in Lemma
\ref{lemma} that the singularities of $X'$ are outside
 the singularities of $\tilde{F}_{3;2}$. Then arguing as in
Theorem \ref{nodal} we prove that they are $30$ ordinary double
points. Let $X''$ be the blow up of $X'$ along $\tilde{F}_{3;2}$.
We deduce as in the proof of Theorem \ref{rank Pic} that
$\rho(X'')=2$.

Analogously we find that after a suitable change of coordinates
the ideal defining $\tilde{F}_{5;1} \subset \mathbb{P}^7$ is
defined by $p^2-wq ,wp-vq ,vp-uq ,up-tq ,w^2-uq ,vw-tq ,uw-tp
,tw-2xp+2zp-tp-2zq ,v^2-tp ,uv-tw ,tv-2xw+2zw-2xp-tp-2zq, u^2-tv
,tu-2xv+2zv-2xw-2yp-tp-2zq ,t^2-2xu+2zu-2xv-2xw-2xp-tp-2zq$ (it
can be shown that $\tilde{F}_{5;1}$ is a degeneration of
$\tilde{F}_{3;2}$, see \cite{Wilson}).
\end{proof}
A natural question arises:
\begin{pro} Is there a primitive contraction of type II with
exceptional divisor isomorphic to $\mathbb{P}^2$ blown up in one
point?
\end{pro}
The answer is probably ``yes'' but an example is difficult to find
because of the high codimension of the anti-canonical image of
this surface.
\begin{rem}
Note that in all the above cases the singular del Pezzo surface
$D'\subset X'$ can be smoothed by a deformation of $X'$. We recall
that Wilson asked in \cite{Wilson} whether this is always possible
for primitive contractions.
\end{rem}

\vskip10pt Department of Mathematics and Informatics, Jagiellonian
University, Reymonta 4, 30-059 Krak\'{o}w, Poland.\\
 \emph{E-mail address:} grzegorz.kapustka@im.uj.edu.pl


\begin{thebibliography}{GPS}
\bibitem[AJ]{AJ} Arapura, D., Jaffe, D.~B., \emph{On Kodaira vanishing for singular varieties}, Proc. of the A.M.S. 105, no. 4 (1989), 911--916..
\bibitem[AR]{AR} Alzati, A., Russo, F., \emph{Some extremal contractions between smooth varieties arising
          from projective geometry}, Proc. London Math. Soc. 89 (2004), 25--53.
\bibitem[Be]{Bertin} Bertin, M.~A., \emph{Examples of Calabi--Yau 3-folds of $\mathbb{P}^7$ with $\rho=1$}, arXiv:math/0701511v1
        [math.AG] (2007).

\bibitem[DH]{DH} Diaz, S., Harbater, D., \emph{Strong Bertini
                theorems}, Trans. Amer. Math. Soc. 324, no. 1 (1991), 73--86.
\bibitem[ES]{EV} van Enckevort, C., van Straten, D., \emph{Monodromy calculations of fourth order equations of Calabi--Yau type},
         Mirror symmetry. V, AMS/IP Stud. Adv. Math., 38, Amer. Math. Soc., Providence, RI (2006), 539--559.

\bibitem[Fr]{Fr} Friedman, R., \emph{Simultaneous Resolution of Threefold
                Double Points}, Math. Ann. 274 (1986), 671--689.
                      \bibitem[Ful]{Fulton} Fulton, W., \emph{Intersection theory}, Ergebnisse der Mathematik und ihrer Grenzgebiete, 2. Springer-Verlag, Berlin, (1984).
\bibitem[GPS]{GPS}
Greuel, G.~M., Pfister G., Sch\"onemann, H.
\newblock {{\sc Singular} 3.0}. A Computer Algebra System for
Polynomial Computations.
\newblock Centre for Computer Algebra, University of Kaiserslautern (2005).
\newblock {\tt http://www.singular.uni-kl.de}.
                \bibitem[Gr1]{Gross1}  Gross, M., \emph{Deforming Calabi--Yau threefolds}, Math. Ann. 308 (1997), 187--220.
\bibitem[Gr2]{Gross2}  Gross, M., \emph{Primitive Calabi--Yau Threefolds},
                J. Differential Geometry 45 (1997), 288--318.
                 \bibitem[GP]{GP} Gross, M., Popescu, S., \emph{Calabi--Yau Threefolds and Moduli of
                Abelian Surfaces I}, Compositio Math. 127 (2001),no. 2, 169--228.
\bibitem[Har]{Har} Hartshorne, R., \emph{Algebraic geometry}, Graduate Texts in Mathematics, no. 52. Springer-Verlag, New York-Heidelberg, (1977).

 \bibitem[KK]{GM} Kapustka, G., Kapustka, M.,
                \emph{Primitive contractions of Calabi--Yau
                threefolds I}, arXiv:math/0703810v3 [math.AG] (2007), to appear in Comm.
                Alg..
 \bibitem[K-P]{KM} Kleppe, J.O., Migliore, J.C, Mir\'{o}-Roig, R., Nagel, U., Peterson, C., \emph{Gorenstein Liaison, Complete Intersection Liaison Invariants and Unobstructedness}, Memoirs AMS, 732, (2001), vol. 154.
 \bibitem[KN]{KN} Kawamata, Y., Namikawa, Y., \emph{Logarithmic deformations of normal crossing varieties and smoothing of degenerate Calabi-Yau varieties},
Invent. Math. 118 (1994), no. 3, 395--409.
\bibitem[Kol]{KOL} Koll\'{a}r, J., \emph{Singularities of
                pairs}, Algebraic geometry---Santa Cruz 1995, Proc. Sympos. Pure Math., 62, Part 1, Amer. Math. Soc., Providence, RI (1997), 221--287.

\bibitem[L1]{Lee2} Lee, N.~H., \emph{Calabi-Yau construction by smoothing normal crossing
                varieties}, arXiv:math/0604596v4 [math.AG] (2007).
\bibitem[L2]{Lee1} Lee, N.~H., \emph{Some Calabi-Yau coverings over singular varieties and new Calabi-Yau threefolds with Picard rank
                one}, arXiv:math/0610060v4 [math.AG] (2008).
\bibitem[LO]{LO} Lee, N.~H., Oguiso, K., \emph{Connecting certain rigid birational non-homeomorphic Calabi--Yau threefolds via
Hilbertscheme}, arXiv:math/0703315v2 [math.AG] (2007).

\bibitem[N1]{namikawa1} Namikawa,~Y., \emph{Deformation theory of Calabi--Yau threefolds and certain invariants of
singularities}, J. Alg. Geometry 6 (1997), 753--776.
\bibitem[N2]{namikawa top} Namikawa,~Y., \emph{Stratified local
moduli of Calabi--Yau
                threefolds}, Topology 41 (2002), 1219-1237.
\bibitem[RS]{RS} Ravindra, G.V., Srinivas, V., \emph{The Grothendieck--Lefschetz
theorem for normal projective varieties}, J. Algebraic Geom. 15
(2006), 563-590.
\bibitem[R1]{Reid}Reid, M., \emph{The moduli space of 3-folds with K = 0 may nevertheless be
                irreducible}, in Hirzebruch Festschrift, Math. Ann. 278 (1987), 329--334.
\bibitem[R2]{R1} Reid, M., \emph{Nonnormal del Pezzo surfaces},  Publ. Res. Inst. Math. Sci. 30 (1994),  no. 5, 695--727.

 \bibitem[R3]{Reid ch}Reid, M., \emph{Chapters on algebraic surfaces}, in Complex algebraic varieties,
                J. Koll\'ar Ed., IAS/Park City lecture notes series (1993 volume), AMS, 1997, 1--154.
\bibitem[R4]{Reid 1} Reid, M., \emph{Constructing algebraic varieties via commutative algebra}, in Proc. of
                4th European Congress of Math
                (Stockholm 2004), European Math Soc. (2005), 655--667.
\bibitem[R5]{Reid 3} Reid, M., \emph{ Minimal models of canonical 3-folds}, Algebraic varieties and
        analytic varieties. litaka, S., ed., Advanced Studies in Pure Math. vol 1.


\bibitem[T]{T} Takagi, H., \emph{On classification of $\mathbb{Q}$-Fano 3-folds of Gorenstein index 2 I, II },  Nagoya Math. J.  167  (2002), 117--155, 157--216.

\bibitem[Weh]{Weh} Wehler, J., \emph{Cyclic Coverings: Deformation and Torelli
       Theorem}, Math. Ann. 274 (1986), 443--472.

 \bibitem[Wi1]{Wilson1} Wilson, P.M.H., \emph{Calabi--Yau manifolds with large Picard number}, Invent. Math. 98 (1989), no. 1, 139--155.
 \bibitem[Wi2]{Wilson2} Wilson, P.M.H., \emph{The K\"{a}hler Cone on Calabi--Yau Threefolds}, Invent.
                Math. 107 (1992), 561--583.
\bibitem[Wi3]{Wilson} Wilson, P.M.H., \emph{Symplectic Deformations of Calabi--Yau
                threefolds},  J. Differential Geom.  45  (1997),  no. 3, 611--637.
\end{thebibliography}
\end{document}